# NON-INTEGRAL TOROIDAL SURGERY ON HYPERBOLIC KNOTS IN $S^3$

CAMERON MCA. GORDON, YING-QING WU, AND XINGRU ZHANG

ABSTRACT. We show that on any hyperbolic knot in $S^3$ there is at most one non-integral Dehn surgery which yields a manifold containing an incompressible torus.

Let $K$ be a knot in the 3-sphere $S^3$ and $M = M_K$ the complement of an open regular neighborhood of $K$ in $S^3$. As usual, the set of slopes on the torus $\partial M$ (i.e. the set of isotopy classes of essential simple loops on $\partial M$) is parameterized by

$$\{m/n : m, n \in \mathbf{Z}, n > 0, (m, n) = 1\} \cup \{1/0\},$$

so that $1/0$ is the meridian slope and $0/1$ is the longitude slope. A slope $m/n$ is called *non-integral* if $n \geq 2$. The manifold obtained by Dehn surgery on $S^3$ along the knot $K$ (equivalently, Dehn filling on $M$ along the torus $\partial M$) with slope $m/n$, is denoted by $M(m/n)$. Now suppose that $K \subset S^3$ is a hyperbolic knot, i.e. the interior of $M$ has a complete hyperbolic metric of finite volume. A basic question in Dehn surgery theory is: when can a surgery on $K$ produce a non-hyperbolic 3-manifold? A special case of this question is: when can a surgery on $K$ produce a toroidal 3-manifold, i.e. a 3-manifold which contains an (embedded) incompressible torus? In [GL1], Gordon and Luecke showed that if $m/n$ is a non-integral slope and $M(m/n)$ contains an incompressible torus, then $n = 2$. In this paper, we show that there is at most one such surgery slope.

**Theorem 1.** *For a hyperbolic knot in $S^3$, there is at most one surgery with non-integral slope producing a manifold containing an incompressible torus.*

There are examples of hyperbolic knots in $S^3$ which admit toroidal surgeries with non-integral slopes [EM]. The best known example is the $(-2, 3, 7)$-pretzel knot, on which the surgery with slope $37/2$ gives a toroidal 3-manifold. It is a conjecture ([Go1], [K, Problem 1.77]) that for a hyperbolic knot in $S^3$, there is at most one surgery with non-integral slope producing a non-hyperbolic 3-manifold, and further, if there is such a surgery, it must be a toroidal surgery and the knot must belong to the collection of examples given in [EM]. Theorem 1 provides some supporting evidence for this conjecture.

1991 *Mathematics Subject Classification.* Mathematics Subject Classification: 57N10, 57M25.
Gordon was Partially supported by an NSF grant. Research at MSRI supported in part by NSF grant #DMS 9022140.

We now go on to the proof of Theorem 1. Our argument is based on applications of results and combinatorial techniques developed in [CGLS,GL1,Go2]. Recall that the *distance* between two slopes $m_1/n_1$ and $m_2/n_2$ is defined as the number $\Delta = \Delta(m_1/n_1, m_2/n_2) = |m_1 n_2 - n_1 m_2|$, which is equal to the minimal geometric intersection number between simple loops representing the two slopes on $\partial M$. Suppose that there are two slopes $m_1/2$ and $m_2/2$ such that both $M(m_1/2)$ and $M(m_2/2)$ contain incompressible tori. It follows from [Go2, Theorem 1.1] that there are four hyperbolic manifolds with toroidal fillings at distance more than 5, but for homological reasons, only one of these is the complement of a knot in $S^3$, namely the figure 8 knot complement, and by [Th] every nonintegral surgery on this manifold is hyperbolic. Hence we have $\Delta(m_1/2, m_2/2) = |2m_1 - 2m_2| = 2|m_1 - m_2| \leq 5$. Note that both $m_1$ and $m_2$ are odd integers. So $|m_1 - m_2|$ is even and thus must be equal to 2. Hence the distance between the two slopes is exactly 4. Our task here is to show that this is impossible. Note, however, that 4 can be realized as the distance between integral toroidal surgery slopes for a hyperbolic knot in $S^3$. For instance the slopes 16 and 20 are both toroidal surgery slopes for the $(-2,3,7)$-pretzel knot. The reason that distance 4 is impossible in our situation is mainly due, as we will see, to the fact that the first homology of $M(m_i/2)$ with $\mathbf{Z}_2$ coefficients is trivial for $i = 1, 2$.

By [GL1] and [GL2] (see [GL1, Theorem 1.2]), for $i = 1, 2$, there is an incompressible torus $\widehat{T}_i$ in $M(m_i/2)$ such that $M \cap \widehat{T}_i = T_i$ is an incompressible, $\partial$-incompressible, twice punctured torus properly embedded in $M$ with each component of $\partial T_i$ having slope $m_i/2$ in $\partial M$. Note that $T_i$ separates $M$ since $M(m_i/2)$ has finite first homology. By an isotopy of $T_i$, we may assume that $T_1$ and $T_2$ intersect transversely, and $T_1 \cap T_2$ has the minimal number of components. So $T_1 \cap T_2$ is a set of finitely many circle components and arc components properly embedded in $T_i$, $i = 1, 2$. Furthermore, no circle component of $T_1 \cap T_2$ bounds a disk in $T_i$ and no arc component of $T_1 \cap T_2$ is boundary parallel in $T_i$, $i = 1, 2$, since $T_i$ is incompressible and $\partial$-incompressible.

We shall use the indices $i$ and $j$ to denote 1 or 2, with the convention that, when they are used together, $\{i, j\} = \{1, 2\}$ as a set.

Let $V_i$ denote the solid torus that is attached to $M$ in forming $M(m_i/2)$. The torus $\widehat{T}_i$ intersects $V_i$ in two disks $B_{i(1)}$ and $B_{i(2)}$, which cut $V_i$ into two 2-handles, which we denote by $H_{i(1)}$ and $H_{i(2)}$. Correspondingly, $\partial T_i = \partial B_{i(1)} \cup \partial B_{i(2)}$ cuts $\partial M$ into two annuli $A_{i(1)}$ and $A_{i(2)}$, where $A_{i(k)} \subset \partial H_{i(k)}$. Each $\partial T_j \cap A_{i(k)}$ consists of exactly 8 essential arcs on $A_{i(k)}$.

The torus $\widehat{T}_i$ separates $M(m_i/2)$ into two submanifolds which we denote by $\widehat{X}_{i(1)}$ and $\widehat{X}_{i(2)}$, with $H_{i(k)} \subset \widehat{X}_{i(k)}$. Correspondingly $T_i$ separates $M$ into two pieces, denoted by $X_{i(1)}$ and $X_{i(2)}$. Thus $X_{i(k)} = M \cap \widehat{X}_{i(k)}$. Note that $\widehat{X}_{i(k)}$ is obtained by attaching the 2-handle $H_{i(k)}$ to $X_{i(k)}$ along the annulus $A_{i(k)}$, $k = 1, 2$. Also note that $F_{i(k)} = \partial X_{i(k)} = T_i \cup A_{i(k)}$ is a closed surface of genus two, $k = 1, 2$.

**Lemma 2.** *For $i = 1, 2$, $k = 1, 2$, we have $H_1(\widehat{X}_{i(k)}, \widehat{T}_i; \mathbf{Z}_2) = 0$.*

*Proof.* Since $H_1(M(m_i/2); \mathbf{Z}_2) = 0$, and $\widehat{T}_i$ is connected, we have

$$0 = H_1(M(m_i/2), \widehat{T}_i; \mathbf{Z}_2) = H_1(\widehat{X}_{i(1)}, \widehat{T}_i; \mathbf{Z}_2) \oplus H_1(\widehat{X}_{i(2)}, \widehat{T}_i; \mathbf{Z}_2),$$

hence each $H_1(\widehat{X}_{i(k)}, \widehat{T}_i; \mathbf{Z}_2) = 0$. □

Now, as in [CGLS, 2.5] and [Go2], we construct two graphs $\Gamma_1$ and $\Gamma_2$ in $\widehat{T}_1$ and $\widehat{T}_2$ respectively by taking the arc components of $T_1 \cap T_2$ as edges and $\widehat{T}_i - int(T_i) = B_{i(1)} \cup B_{i(2)}$ as (fat) vertices. So $\Gamma_i$ is a graph on a torus with two vertices. The exterior of the graph $\Gamma_i$ in $\widehat{T}_i$ is the set of faces of $\Gamma_i$. Each face of the graph $\Gamma_i$ is a surface properly embedded in $X_{j(k)}$ for some $k = 1$ or 2. Since each component of $\partial T_1$ intersects each component of $\partial T_2$ in exactly 4 points, each of the two vertices of $\Gamma_i$ has valency (i.e. the number of edge endpoints at the vertex) 8 and there are exactly 8 edges in $\Gamma_i$. If $e$ is an edge with an endpoint at a vertex $B_{i(k)}$ of $\Gamma_i$ ($k = 1$ or 2), then that endpoint is in $\partial B_{i(k)} \cap \partial B_{j(l)}$ for some vertex $B_{j(l)}$ ($l = 1$ or 2) of $\Gamma_j$, and the endpoint of $e$ at the vertex $B_{i(k)}$ is given the label $l$. So the labels of edge endpoints around each of the vertices $B_{i(k)}$, $i = 1, 2$, $k = 1, 2$, are as shown in Figure 1.

Figure 1

A *cycle* in $\Gamma_i$ is a subgraph homeomorphic to a circle (where the fat vertices of $\Gamma_i$ are considered as points). The number of edges in a cycle is its *length*. A length one cycle is also called a *loop*. A loop of a graph is said to be *trivial* if it bounds a disk face of the graph. So $\Gamma_i$ has no trivial loops. We shall also consider two parallel loops as forming a cycle of length two. (Two edges of $\Gamma_i$ are said to be *parallel* in $\Gamma_i$ if they, together with two arcs in $\partial T_i$, bound a disk in $T_i$.) A cycle in $\Gamma_i$ consisting of two parallel adjacent loops is called an *S-cycle*. This definition is a specialization to our current situation of the usual definition of an S-cycle given in [W, GL1]. The disk face of $\Gamma_i$ bounded by an S-cycle in $\Gamma_i$ is called the *S-disk* of the S-cycle.

Figure 2

Let $\overline{\Gamma}_i$ be the reduced graph of $\Gamma_i$, i.e. the graph obtained from $\Gamma_i$ by amalgamating each complete set of mutually parallel edges of $\Gamma_i$ to a single edge. Then up to homeomorphism of $\widehat{T}_i$, $\overline{\Gamma}_i$ is a subgraph of the graph illustrated in Figure 2 (for a proof of this, see [Go2, Lemma 5.2]). It follows in particular that in $\Gamma_i$ the number of loops at the vertex $B_{i(1)}$ is equal to the number of loops at the vertex $B_{i(2)}$, and thus the number of loops of $\Gamma_i$ is even. From now on we will take Figure 2 as a fixed model graph of which our $\overline{\Gamma}_i$ is a subgraph. Note that each non-loop edge in $\overline{\Gamma}_i$ represents at most two edges of $\Gamma_i$, otherwise there would be a pair of edges parallel in both $\Gamma_i$, contradicting Lemma 3(1).

If $\Gamma_i$ contains some loops, then they cut $T_i$ into disks together with two annuli, which will be called the *loop complement annuli*.

**Lemma 3.** (1) *No pair of edges can be parallel in both $\Gamma_1$ and $\Gamma_2$.*

(2) *Any closed curve $\alpha$ on $T_i$ intersects $\Gamma_i$ in an even number of points. In particular, any loop complement annulus contains an even number of (non-loop) edges.*

(3) (The parity rule) *A component $e$ of $T_1 \cap T_2$ is a loop on $\Gamma_1$ if and only if it is a non-loop on $\Gamma_2$. An edge of $\Gamma_i$ has the same label at its two endpoints if and only if it is not a loop.*

*Proof.* (1) Otherwise, the manifold $M$ would be cabled by [Go2, Lemma 2.1], contradicting the assumption that $M$ is hyperbolic.

(2) We have seen above that $\widehat{T}_j$ is a separating torus, hence $\alpha$ intersects $\widehat{T}_j$ in an even number of points. Since $\alpha \cap \Gamma_i = \alpha \cap \widehat{T}_j$, the result follows. For the statement about a loop complement annulus $A$, just take $\alpha$ to be the central curve of $A$.

(3) We fix an orientation on $T_i$ and let each component of $\partial T_i$ have the induced orientation. Then the two components of $\partial T_i$ are not homologous in $\partial M$ since $T_i$ separates $M$. (The corresponding vertices of $\Gamma_i$ are *anti-parallel* in the terminology of [CGLS, p. 278]). Since $M$, $T_1$ and $T_2$ are all orientable, the parity rule given in [CGLS, Page 279] holds and has the above special form in our present situation. □

**Lemma 4.** $\Gamma_i$ *cannot have two disk faces $D_1, D_2$, each with an even number of edges, such that $\partial D_1$ and $\partial D_2$ are nonparallel, nonseparating curves on $F_{j(k)}$ for some $k = 1, 2$.*

*Proof.* Otherwise, since $F_{j(k)}$ has genus 2, the disks $D_i$ cuts $X_{j(k)}$ into a manifold with sphere boundary, which must be a 3-ball $B$ because $X_{j(k)}$, as a submanifold of $S^3$ with a single boundary component, is irreducible. Thus we have

$$\widehat{X}_{j(k)} = \widehat{T}_j \cup H_{j(k)} \cup D_1 \cup D_2 \cup B.$$

Now $H_1(\widehat{T}_j \cup H_{j(k)}, \widehat{T}_j; \mathbf{Z}_2) = \mathbf{Z}_2$ is generated by the core of $H_{j(k)}$. Since $\partial D_i$ has an even number of edges, it runs over $H_{j(k)}$ an even number of times, so it represents zero in $H_1(\widehat{T}_j \cup H_{j(k)}, \widehat{T}_j; \mathbf{Z}_2)$. Therefore,

$$H_1(\widehat{X}_{j(k)}, \widehat{T}_{j(k)}; \mathbf{Z}_2) = H_1(\widehat{T}_j \cup H_{j(k)}, \widehat{T}_{j(k)}; \mathbf{Z}_2)/\langle \partial D_1, \partial D_2 \rangle = \mathbf{Z}_2/\langle 0, 0 \rangle = \mathbf{Z}_2.$$

This contradicts Lemma 2, completing the proof. □

**Lemma 5.** *Let $\partial_i$ be a boundary component of $T_i$, $i = 1, 2$. If the four points of $\partial_1 \cap \partial_2$ appear in the order $x_1, x_2, x_3, x_4$ on $\partial_1$, then they also appear in the same order on $\partial_2$, in some direction. In particular, if two of the points $x_p, x_q$ are adjacent on $\partial_1$ among the four points, then they are also adjacent on $\partial_2$.*

*Proof.* We may choose coordinates on $\partial M$ so that $\partial_1$ is the $1/0$ curve and $\partial_2$ is the $1/4$ curve. The lemma is obvious by drawing such curves on a torus. □

We will often apply this lemma to the endpoints of a pair of edges $e_1, e_2$. If each $e_i$ has an endpoint $x_i$ at the vertex $B_{1(l)}$, with label $k$, then at $B_{2(k)}$ $e_1$ and $e_2$ both have label $l$. The above says that $x_1, x_2$ are adjacent on $\partial B_{1(l)}$ among all edge endpoints labeled $k$ if and only if they are adjacent on $\partial B_{2(k)}$ among all edge endpoints labeled $l$.

By the parity rule, one of $\Gamma_1$ and $\Gamma_2$, say $\Gamma_1$, contains at least 4 loops. So $\Gamma_1$ has either 8, 6 or 4 loops.

**Lemma 6.** $\Gamma_i$ *cannot have 8 loops.*

*Proof.* Suppose $\Gamma_1$, say, has 8 loops. Then the four loops based at the vertex $B_{1(1)}$ in $\Gamma_1$ form three S-cycles, bounding three S-disks $D_1, D_2$ and $D_3$. The disks $D_1$ and $D_3$ are on the same side of $T_2$, say $X_{2(1)}$. The four corresponding edges of $\Gamma_2$ connect the two vertices $B_{2(1)}$ and $B_{2(2)}$, and by Lemma 3(1) they are mutually nonparallel in $\Gamma_2$. Thus $\partial D_1$ and $\partial D_3$ are nonparallel curves on the genus 2 surface $F_{2(1)} = \partial X_{2(1)}$. Since each $\partial D_i$ is an S-cycle, it runs over $H_{2(1)}$ twice in the same direction, so it is a nonseparating curve on $F_{2(1)}$. This contradicts Lemma 4. □

**Lemma 7.** $\Gamma_i$ *cannot have 6 loops.*

*Proof.* Suppose $\Gamma_1$, say, has 6 loops. By Lemma 3(2) the number of non-loop edges in each loop complement annulus $A_k$ is even. Thus the two non-loop edges must be on the same $A_k$. Up to isomorphism, the graph $\Gamma_1$ is one of the two graphs shown in Figure 3(a)–(b). However, Figure 3(a) is not possible because the longitude of the torus intersects the graph in an odd number of points, which contradicts Lemma 3(2).

Now consider $\Gamma_2$. For the same reason as above, the number of non-loop edges on each loop complement annulus $A_k$ is even, so we have two possibilities, shown in Figure 3(c)–(d). Again, Figure 3(c) can be ruled out by looking at the intersection number between $\Gamma_2$ and a longitude of the torus.

Consider the two edges $e_1, e_2$ in Figure 3(d). The two endpoints $x_1, x_2$ of $e_1, e_2$ on $\partial B_{2(2)}$ are adjacent among all points of $\partial B_{1(1)} \cap \partial B_{2(2)}$, so by Lemma 5 they are also adjacent on $\partial B_{1(1)}$. From Figure 3(b) we see that the other two endpoints $y_1, y_2$ of $e_1, e_2$ are also adjacent on $\partial B_{1(1)}$ among all points of $\partial B_{1(1)} \cap \partial B_{2(1)}$. However, on Figure 3(d) the two endpoints of $e_1, e_2$ on $\partial B_{2(1)}$ are not adjacent among all edge endpoints labeled 1, which contradicts Lemma 5. □

Figure 3

We may now assume that $\Gamma_1$ has exactly four loops. By the parity rule $\Gamma_2$ also has exactly four loops. By Lemma 3(2) each loop complement annulus $A_k$ contains an even number of edges, so there are four possible configurations for the graph $\Gamma_i$, according to whether each $A_k$ contains exactly two edges, and if there are two, whether they are parallel. See Figure 4(a)–(d). However, Figure 4(a) is impossible for each of $\Gamma_1$ and $\Gamma_2$ because the longitude of the torus intersects the graph in an odd number of points.

Figure 4

**Lemma 8.** *Figure 4(b) is impossible for each of $\Gamma_1$ and $\Gamma_2$.*

*Proof.* Suppose that one of $\Gamma_1$ and $\Gamma_2$, say $\Gamma_1$, is as shown in Figure 4(b). Let $D_1$ be the $S$-disk bounded by the $S$-cycle $\{e_1, e_2\}$, and $D_2$ the disk bounded by the cycle $\{e_3, e_4\}$, as shown in Figure 5.

Figure 5

Then $D_1$ and $D_2$ are contained in the same side of $T_2$, say in $X_{2(1)}$. In $\Gamma_2$, the two edges $\{e_1, e_2\}$ also form a cycle which we denote by $\sigma_1$. There are two possibilities for $\sigma_1 = \{e_1, e_2\}$, depending on whether or not $e_1$ and $e_2$ lie on the same loop complement annulus. See Figure 6.

Figure 6

It is clear that $\partial D_1$ and $\partial D_2$ are nonparallel curves on $F_{2(1)}$. Also, as the boundary of an S-disk, $\partial D_1$ is always nonseparating. Now $\partial D_2$ consists of one loop at each vertex of $\Gamma_2$, together with two arcs on the boundary of the handle $H_{2(1)}$. In the case of Figure 6(a), there is an arc with its endpoints on $\partial D_1$, intersecting $\partial D_2$ at a single point, and hence $\partial D_2$ is also nonseparating. This contradicts Lemma 4.

In the case of Figure 6(b), the two endpoints of $e_1$ and $e_2$ on $\partial B_{2(1)}$ are labeled 1, say, and are adjacent among all points of $\partial B_{1(1)} \cap \partial B_{2(1)}$ on $\partial B_{2(1)}$. But from Figure 5 we can see that they are not adjacent among all edge endpoints labeled 1 on $\partial B_{1(1)}$, which contradicts Lemma 5.  □

**Lemma 9.** $\Gamma_i$ *cannot be as in Figure 4(c).*

*Proof.* Suppose $\Gamma_1$, say, is as shown in Figure 4(c). Consider the disks $D_1, D_2$ shown in Figure 7. They are on the same side of $T_2$, say in $X_{2(1)}$. The curve $\partial D_1$ is an S-cycle on $\Gamma_1$, so it is automatically nonseparating on $F_{2(1)}$. The curve $\partial D_2$ is also nonseparating, because it runs over $H_{2(1)}$ three times. Clearly they are nonparallel, so by the proof of Lemma 4 we see that

$$\widehat{X}_{2(1)} = \widehat{T}_2 \cup H_{2(1)} \cup D_1 \cup D_2 \cup B,$$

where $B$ is a 3-ball. Similarly

$$\widehat{X}_{2(2)} = \widehat{T}_2 \cup H_{2(2)} \cup D_3 \cup D_4 \cup B.$$

Hence

$$H_1(M(m_2/2); \mathbf{Z}_2) = H_1(\widehat{T}_2 \cup V_2; \mathbf{Z}_2) / \langle \partial D_1, \partial D_2, \partial D_3, \partial D_4 \rangle = (\mathbf{Z}_2)^4 / \langle \partial D_1, \ldots \partial D_4 \rangle.$$

Each of $\partial D_1$ and $\partial D_3$ runs over $H_{2(k)}$ an even number of times, so if we denote by $\sigma_1, \sigma_3$ the corresponding cycles on $\Gamma_2$ with each fat vertex shrunk to a point, then homologically we have $\partial D_i = \sigma_i$ in $H_1(\widehat{T}_2 \cup V_2; \mathbf{Z}_2)$, for $i = 1, 3$.

Figure 7

Consider $\Gamma_2$. We have shown so far that $\Gamma_2$ must be the graph in either Figure 4(c) or 4(d). In either case we have $\sigma_1 + \sigma_3 = 0$ in $H_1(\widehat{T}_2; \mathbf{Z}_2)$, so homologically we have $\partial D_1 = \partial D_3$. Therefore,

$$0 = H_1(M(m_2/2); \mathbf{Z}_2) = (\mathbf{Z}_2)^4 / \langle \partial D_1, \partial D_2, \partial D_4 \rangle \neq 0.$$

This contradiction completes the proof of the lemma. □

**Lemma 10.** *$\Gamma_1$ and $\Gamma_2$ cannot both be as in Figure 4(d).*

*Proof.* In this case the graph $\Gamma_i$ is contained in an annulus $Q_i$ in $\widehat{T}_i$, so we can redraw the graph on $Q_i$, $i = 1, 2$, as in Figure 8. Orient the non-loop edges in both $\Gamma_i$ so that they go from $B_{i(1)}$ to $B_{i(2)}$. This determines the orientations of the loop edges:

they must go from the endpoints labeled 1 to those labeled 2. Use $e_i^-, e_i^+$ to denote the tail and head of $e_i$, respectively.

Label the edges of $\Gamma_1$ as in Figure 8(a). There is a correspondence between the edges of $\Gamma_1$ and $\Gamma_2$. Up to isomorphism we may assume that $e_1$ on $\Gamma_2$ is as shown in Figure 8(b). There are only two non-loop edges on $\Gamma_2$ labeled 2. They correspond to loops in $\Gamma_1$ based at $B_{1(2)}$, so the other one must be $e_2$. Now consider the four points of $\partial B_{1(2)} \cap \partial B_{2(1)}$. On $\partial B_{1(2)}$ they appear in the order $e_1^-, e_2^-, e_7^+, e_3^+$, so by Lemma 5 they appear in the same order on $\partial B_{2(1)}$. This determines the edges $e_7$ and $e_3$. Similarly one can determine the correspondence between the other edges. The labeling of the edges are shown in Figure 8.

Figure 8

Let $D_1, D_2$ be the two disk faces of $\Gamma_1$ shown in Figure 8(a). Note that they are on different sides of $T_2$. Since $D_1$ is an S-disk, and $D_2$ has an odd number of edges, both $\partial D_1$ and $\partial D_2$ are nonseparating curves on $F_{2(1)}$ and $F_{2(2)}$, respectively. Let $P$ be the twice punctured annulus obtained from $Q_2$ by removing the interiors of the fat vertices. Consider the complex $Y = \partial M \cup P \cup D_1 \cup D_2$. Denote by $N(S)$ the regular neighborhood of a set $S$ in $M$. Notice that $N(\partial M \cup P)$ has boundary the torus $\partial M$ and a surface $F$ of genus three. We have shown that $\partial D_1, \partial D_2$ are nonseparating. Since they lie on different sides of $P$, they are also nonparallel. Thus after adding the two 2-handles $N(D_1)$ and $N(D_2)$ to $N(\partial M \cup P)$, the manifold $N(Y)$ has boundary consisting of two tori.

We would like to calculate the homology of $N(Y)$. Consider the complex $P \cup \partial M$ shown in Figure 9(a). We have determined the graph $\Gamma_2$ on $P$. Up to homeomorphism there is a unique way of connecting the endpoints of $e_1, e_2$ by arcs on $A_{2(1)}$ to form the loop $\partial D_1$. Similarly for $\partial D_2$. See Figure 9(a)–(b). This also determines the boundary slope $m_1/2$ on $\partial M$, as shown in Figure 9(c). Now pick a basis $\{u, v, x_1, x_2\}$ for $H_1(\partial M \cup P) = \mathbf{Z}^4$ as in Figure 9(d), where $u$ is the loop $\partial B_{2(1)}$, $v$ is the inner boundary circle of $Q_2$, and $x_j$ runs over the annulus $A_{2(j)}$.

Figure 9

Consider the curves $\partial D_1, \partial D_2$ shown in Figure 9(a)–(b). Calculating their homology classes in $H_1(\partial M \cup P)$, we have

$$\begin{aligned} \partial D_1 &= x_1 + v + x_1 - u = 2x_1 + v - u, \\ \partial D_2 &= -x_2 + v + x_2 + u + v - x_2 - u = 2v - x_2. \end{aligned}$$

The boundary slope $m_2/2$ is represented by $u$. From Figure 9(c) we can see that the slope $m_1/2$ is a $(1, -4)$ curve, represented by $u - 4(x_1 + x_2)$. The meridian slope $r_0 = 1/0$ of $\partial M$ is characterized by the property that it has geometric intersection number 2 with each of $m_1/2$ and $m_2/2$, so it must be a $(1, -2)$ curve, represented by $u - 2(x_1 + x_2)$. Thus in homology we have

$$r_0 = u - 2(x_1 + x_2).$$

Denote by $W$ the manifold obtained from $N(Y)$ by Dehn filling on $\partial M$ along the meridian slope. That is, $W = N(Y)(1/0) = N(Y)(r_0)$. Then $W$ is obtained from $\partial M \cup P$ by adding two 2-handles along $\partial D_1, \partial D_2$, then adding a solid torus along the slope $r_0$. So $H_1(W)$ has the presentation

$$\begin{aligned} \langle u, v, x_1, x_2 \,:\, 2x_1 - u + v = 2v - x_2 &= u - 2(x_1 + x_2) = 0\rangle \\ = \langle x_1, x_2 \,:\, 3x_2 = 0\rangle &= \mathbf{Z} \oplus \mathbf{Z}_3. \end{aligned}$$

We have seen that $\partial N(Y)$ consists of two tori, hence $W = N(Y)(1/0)$ has boundary a single torus. On the other hand, $W$ is a submanifold of $M(r_0) = S^3$ with boundary a torus, hence $H_1(W) = \mathbf{Z}$. This contradiction completes the proof of the lemma, hence the proof of Theorem 1. □

CAMERON MCA. GORDON, DEPARTMENT OF MATHEMATICS, UNIVERSITY OF TEXAS AT AUSTIN, AUSTIN, TX 78712
   *Current address*: MSRI, 1000 Centennial Drive, Berkeley, CA 94720-5070
   *E-mail address*: gordon@math.utexas.edu

YING-QING WU, DEPARTMENT OF MATHEMATICS, UNIVERSITY OF IOWA, IOWA CITY, IA 52242
   *Current address*: MSRI, 1000 Centennial Drive, Berkeley, CA 94720-5070
   *E-mail address*: wu@math.uiowa.edu

XINGRU ZHANG, DEPARTMENT OF MATHEMATICS, OKLAHOMA STATE UNIVERSITY, STILLWATER, OK 74078
   *E-mail address*: xingru@math.okstate.edu